      \documentstyle[12pt,graphicx]{article}
\begin{document}
\baselineskip 18pt
\def\today{\ifcase\month\or
 January\or February\or March\or April\or May\or June\or
 July\or August\or September\or October\or November\or December\fi
 \space\number\day, \number\year}
\def\thebibliography#1{\section*{References\markboth
 {References}{References}}\list
 {[\arabic{enumi}]}{\settowidth\labelwidth{[#1]}
 \leftmargin\labelwidth
 \advance\leftmargin\labelsep
 \usecounter{enumi}}
 \def\newblock{\hskip .11em plus .33em minus .07em}
 \sloppy
 \sfcode`\.=1000\relax}
\let\endthebibliography=\endlist
\def\lsim{\ ^<\llap{$_\sim$}\ }
\def\gsim{\ ^>\llap{$_\sim$}\ }
\def\r2{\sqrt 2}
\def\beq{\begin{equation}}
\def\eeq{\end{equation}}
\def\beqn{\begin{eqnarray}}
\def\eeqn{\end{eqnarray}}
\def\rmuu{\gamma^{\mu}}
\def\rmud{\gamma_{\mu}}
\def\PL{{1-\gamma_5\over 2}}
\def\PR{{1+\gamma_5\over 2}}
\def\sinW2{\sin^2\theta_W}
\def\AEM{\theta_{_{n}}_{EM}}
\def\mul{M_{\tilde{u} L}^2}
\def\mur{M_{\tilde{u} R}^2}
\def\mdl{M_{\tilde{d} L}^2}
\def\mdr{M_{\tilde{d} R}^2}
\def\mz2{M_{z}^2}
\def\c2b{\cos 2\beta}
\def\au{A_u}
\def\ad{A_d}
\def\cob{\cot \beta}
\def\v#1{v_#1}
\def\tb{\tan\beta}
\def\epem{$e^+e^-$}
\def\KK{$K^0$-$\bar{K^0}$}
\def\wi{\omega_i}
\def\xj{\chi_j}
\def\Wmu{W_\mu}
\def\Wnu{W_\nu}
\def\m#1{{\tilde m}_#1}
\def\mH{m_H}
\def\mw#1{{\tilde m}_{\omega #1}}
\def\mx#1{{\tilde m}_{\chi^{0}_#1}}
\def\mc#1{{\tilde m}_{\chi^{+}_#1}}
\def\mwi{{\tilde m}_{\omega i}}
\def\mxi{{\tilde m}_{\chi^{0}_i}}
\def\mci{{\tilde m}_{\chi^{+}_i}}
\def\mz{M_z}
\def\sw{\sin\theta_W}
\def\cw{\cos\theta_W}
\def\cb{\cos\beta}
\def\sb{\sin\beta}
\def\rwi{r_{\omega i}}
\def\rxj{r_{\chi j}}
\def\rfp{r_f'}
\def\Kik{K_{ik}}
\def\Fq2{F_{2}(q^2)}
\def\f{\({\cal F}\)}
\def\d1{{\f(\tilde c;\tilde s;\tilde W)+ \f(\tilde c;\tilde \mu;\tilde W)}}
\def\tw{\tan\theta_W}
\def\sec2w{sec^2\theta_W}

\begin{titlepage}
\begin{flushright}
\date{\today}
\end{flushright}
\begin{center}
{\LARGE {\bf AREAL OPTIMIZATION OF POLYGONS}}\footnote{submitted for publication}\\
\vskip 0.5 true cm \vspace{2cm}
\renewcommand{\thefootnote}
{\fnsymbol{footnote}}
 \large{\textsl{\textsf{Erica Walker}}$^{a}$, \textsl{\textsf{Raza M.
 Syed}}$^{a,b}$, and \textsl{\textsf{Achille Corsetti}}$^{a,b}$}
\vskip 0.5 true cm
\end{center}

\noindent{$a$. Health Careers Academy, 110 The Fenway, Cahners Hall, Boston, MA 02115\\
{ $b$. Department of Physics, Northeastern University,
360 Huntington Ave., Boston\\$~~~~~$MA 02115-5000}} \\

\vskip 1.0 true cm
\centerline{\bf Abstract}
\medskip
\noindent
 We will first solve the following problem analytically:
given a piece of wire of specified length, we will find where the
wire should be cut and bent to form two regular polygons not
necessarily having the same number of sides, so that the combined
area of the polygons thus formed is maximized, minimized, greater
than, and less than a specified area. We will extend the results
to the cases where the wire is divided into three and finally into
an arbitrary number of segments. The second problem we will solve
is as follows: two wires of specified length are to be bent into
two regular polygons whose total number of sides is fixed.  We
will determine how the total number of polygonal sides are to be
allocated between the wires so that the total area of the polygons
is maximized. We will extend the results found here to the case
where we are given any number of wires of specified length.

\end{titlepage}

\section{Introduction}
A segment of length, $x$, is to be cut from a wire of length, $L$.
This segment, along with the remaining segment of length $L-x$,
are then bent to form two regular polygons not necessarily having
the same number of sides. We will begin by finding the total area,
$A_{_{total}}$, of these two polygons in terms of $x$ and $L$.
From this general form, we will solve the extrema problem; that
is, we will compute $x$-- the place where the wire must be cut and
divided into the perimeter of the polygons so that the total area
of the polygons are maximized and minimized.  Next, we will solve
the inequality problems by again computing $x$, the place where
the wire should be severed so that the total area of these two
polygons is greater than or less than a specified area, ${\mathsf
A}$. In addition to the above 1-cut (2-partitioned) problem, we
will also analyze the above four cases for the 2-cut problem and
finally generalizing it to a ${\mathsf k}$-cut (${\mathsf
k+1}$-partitioned problem).  We will illustrate the above four
cases by providing examples of specific polygons.

    In our second problem, we are given two wires of specified
lengths which we will bend to form two regular polygons with the
total number of polygonal sides fixed.  We will then calculate
ways in which the total number of sides can be allocated between
these two wires so that the combined area of the polygons is
maximized. This will generate a non-linear coupled equation in the
number of polygonal sides. This equation will be solved
numerically using a computer program.  Finally, we will generalize
this problem to ${\mathrm k}$-wires of specified lengths.

     The outline of the paper is as follows: in section 2 we will
examine the 2-partitioned case.  Here we will devote sub-sections
to deriving a general expression for $A_{_{total}}$ in terms of
$x$ and $ L$, computing the exact form for $x$ for which
$A_{_{total}}$ is  maximized and minimized, and considering lower
and upper bound problems as mentioned above. Here, we also mention
constraints on ${\mathsf A}$ and $L$. Finally, we illustrate each
of the above scenarios with specific examples. In sections 3 and 4
we will consider an extension of section 2 examining the
3-partition and the ${\mathsf (k+1)}$-partition cases,
respectively. In section 5 we consider the problem where two wires
of specified length are bent into two regular polygons where the
total number of the polygonal sides is fixed.  Here, we find how
the total number of sides must be divided between the polygons so
that the combined area of the two polygons is maximized.  Next, we
generalize these results to the case of ${\mathrm k}$ wires of
specified lengths. Section 6 is devoted to the conclusion. In
appendix A we show that the maximum area of a polygon for a fixed
perimeter occurs in the limiting case of a circle.  Appendix B
contains the computer program which looks for the right
combination of polygonal sides that would maximize the total area.
.

\section{2-partitioned wire}
{\bf{\textsf{I.~~The total area of two polygons}}}\\
In this sub-section we will derive $A_{_{total}}$, the total
combined area of two regular $m$- and $n$-sided polygons. We will
start by finding the area of each polygon using the formula:
$Area_{_{polygon}} =\frac{1}{2}ap$, where $a$ and $p$ are the
apothem and perimeter of the polygon, respectively. This is
equivalent to finding the area of one triangle (see accompanying
figure) and multiplying it by the number of sides since the number
of congruent triangles is equal to the number of sides.

Cutting length $x$ from a total length, $L$, we form a regular
$n$-gon with side length $\frac{x}{n}$. Dropping a perpendicular
of length $a$ from the center of this $n$-gon to a side will
bisect the side. Further, the radius of the polygon will bisect
the interior angle which is given by $\theta_{_{n}} =
(\frac{1}{2}- \frac{1}{n})\pi$. From the right triangle AOB we
find $a = \frac{x}{2n}\tan{\theta_{_{n}}}$. Hence, the area of the
$n$-sided polygon is given by
$A_{_{n}}(x)=\frac{x^2}{4n}\tan{\theta_{_{n}}}$. Since the
remaining length is $y\equiv L-x$, the area of the $m$-sided
polygon will be given by $A_{_{m}}(y)= \frac{y^2}
{4m}\tan{\theta_{_{m}}}$ where $\theta_{_{m}}
=(\frac{1}{2}-\frac{1}{m})\pi$.  Finally, the total area of the
two polygons is given by the expression:
\begin{equation}
A_{_{total}}(x,y) =
A_{_{n}}(x)+A_{_{m}}(y)=\frac{x^2}{4n}\tan{\theta_{_{n}}} +
\frac{(L-x)^2}{4m}\tan{\theta_{_{m}}}
\end{equation}

In the case of identical polygons ($m=n$ and therefore
$\theta_{_{m}}=\theta_{_{n}}$), Eq.1 reduces to
\begin{equation}
A_{_{total}}^{^{(m=n)}} =
\frac{2x^2-2xL+L^2}{4n}\tan{\theta_{_{n}}}
\end{equation}
{\bf\textsf{II.~~Minimum total area of two polygons}}\\
In this sub-section we consider the minima problem. We begin by
finding the critical values of $A_{_{total}}$. From Eq.(1) we
obtain,
\begin{equation}
\frac{dA_{_{total}}}{dx} = \frac{x\tan{\theta_{_{n}}}}{2n}-
\frac{(L-x)\tan{\theta_{_{m}}}}{2m}=0
\end{equation}
which results in a minimum point:
\begin{equation}
 x_{_{min}} = \frac{nL\tan{\theta_{_{m}}}}{m\tan{\theta_{_{n}}}+ n\tan{\theta_{_{m}}}}
\end{equation}
because,
\begin{equation}
\frac{d^2A_{_{total}}}{dx^2}=
\frac{\tan{\theta_{_{n}}}}{2n}+\frac{\tan{\theta_{_{m}}}}{2m} >0
\end{equation}
This is true since,
$\pi/6\leq\theta_{_{n}},\theta_{_{m}}\leq\pi/2$. Further,
\begin{equation}
 y_{_{min}} = \frac{mL\tan{\theta_{_{n}}}}{m\tan{\theta_{_{n}}}+ n\tan{\theta_{_{m}}}}
\end{equation}
With this result of $x_{_{min}}$ and $ y_{_{min}}$, we now find
expressions for $A_{_{n}}(x_{_{min}})$, $A_{_{m}}(y_{_{min}})$,
and  the total minimum area of the two polygons,
$A_{_{min-total}}$:
\begin{eqnarray}
A_{_{n}}(x_{_{min}})=\frac{L^2n\tan{\theta_{_{n}}}\tan^2{\theta_{_{m}}}}{4(m\tan{\theta_{_{n}}}+
n\tan{\theta_{_{m}}})^2}\nonumber\\
A_{_{m}}(y_{_{min}})=\frac{L^2m\tan^2{\theta_{_{n}}}\tan{\theta_{_{m}}}}{4(m\tan{\theta_{_{n}}}+
n\tan{\theta_{_{m}}})^2}\nonumber\\
A_{_{min-total}} \equiv A_{_{total}}(x_{_{min}},y_{_{min}})=
\frac{L^2\tan{\theta_{_{n}}}\tan{\theta_{_{m}}}}{4(m\tan{\theta_{_{n}}}+
n\tan{\theta_{_{m}}})}
\end{eqnarray}
 Further, in the case when $m=n$, we get:
\begin{eqnarray}
x^{^{(m=n)}}_{_{min}}=y^{^{(m=n)}}_{_{min}}=\frac{L}{2}\nonumber\\
A_{_{n}}(x_{_{min}})=A_{_{m}}(y_{_{min}})=\frac{L^2}{16n}\tan{\theta_{_{n}}}\nonumber\\
 A_{_{min-total}}^{^{(m=n)}}= \frac{L^2}{8n}\tan{\theta_{_{n}}}
\end{eqnarray}
{\bf\textsf{III.~~Maximum total area of two polygons}}\\
This sub-section is devoted to the corresponding maxima problem.
The feasible domain of $x$ is [$0$,$L$] and the only critical
value is $x_{_{min}}$. Therefore, using
\begin{eqnarray}
A_{_{total}}(0,L)=\frac{L^2}{4m}\tan{\theta_{_{m}}}\nonumber\\
A_{_{total}}(x_{_{min}},y_{_{min}})=\frac{L^2}{4m}\tan{\theta_{_{m}}}\left[
\frac{1}{1+\frac{n}{m}\frac{\tan{\theta_{_{m}}}}{\tan{\theta_{_{n}}}}}\right]\nonumber\\
A_{_{total}}(L,0)=\frac{L^2}{4n}\tan{\theta_{_{n}}}
\end{eqnarray}
we can conclude that maximum area occurs when $x=0$ or $x=L$
(endpoint maxima), depending on the integers $m$ and $n$. This is
equivalent to saying all the wire is used for the $m$-sided
polygon (if $m>n$) or $n$-sided polygon (if $n>m$). Thus,
\begin{eqnarray}
A_{_{max-total}}^{^{(m>n)}}=\frac{L^2}{4m}\tan{\theta_{_{m}}}\nonumber\\
A_{_{max-total}}^{^{(n>m)}}=\frac{L^2}{4n}\tan{\theta_{_{n}}}
\end{eqnarray}
{\bf\textsf{IV.~~Lower bound on the area of two polygons}}\\
In this sub-section we will use Eq.(1) to find the possible range
of values of $x$ for which the combined area of the two polygons
exceeds a given area, ${\mathsf A}$. This translates to the
following problem:
\begin{equation}
 \frac{x^2}{4n}\tan{\theta_{_{n}}} +
\frac{(L-x)^2}{4m}\tan{\theta_{_{m}}}>{\mathsf A}
\end{equation}
 which simplifies to
 \begin{equation}
\left(x-\left[x_{_{min}}-\hat{x}\right]\right)\left(x-\left[x_{_{min}}+\hat{x}\right]\right)>0
\end{equation}
 where we have defined for brevity
 \begin{equation}
\hat{x}=\sqrt{\frac{mx_{_{min}}}{nL\tan{\theta_{_{m}}}}\left(4n{\mathsf
A}-Lx_{_{min}}\tan{\theta_{_{n}}}\right)}
\end{equation}
and $x_{_{min}}$ is defined as before. Therefore,
\begin{equation}
x\in
\left(0,x_{_{min}}-\hat{x}\right)\cup\left(x_{_{min}}+\hat{x},L\right)
\end{equation}
\noindent
{\emph{{\textbf{Constraints:}}}\\
Firstly, $\hat{x}$ must be real. This implies
\begin{equation}
L^2-4{\mathsf
A}\left(\frac{n}{\tan{\theta_{_{n}}}}+\frac{m}{\tan{\theta_{_{m}}}}
\right)\leq 0
\end{equation}
Thus,
\begin{equation}
0<L\leq L_1
\end{equation}
where $L_1$ is defined by
\begin{equation}
L_1=2\sqrt{{\mathsf
A}\left(\frac{n}{\tan{\theta_{_{n}}}}+\frac{m}{\tan{\theta_{_{m}}}}
\right)}
\end{equation}
Secondly, $x$ must be non-negative, i.e.,
\begin{equation}
L^2-\frac{4m{\mathsf A}}{\tan{\theta_{_{m}}}}\geq 0
\end{equation}
The solution to this inequality is
\begin{equation}
L\geq L_2
\end{equation}
with $L_2$ given by
\begin{equation}
L_2=2\sqrt{\frac{m{\mathsf A}}{\tan{\theta_{_{m}}}}}
\end{equation}
Finally, noting that $L_1>L_2$ and taking the intersection of the
inequalities \\ (16) and (19) we find that $L$ must satisfy,
\begin{equation}
L_2\leq L \leq L_1
\end{equation}
This is equivalent to
\begin{equation}
{\mathsf A}_2\leq {\mathsf A}\leq {\mathsf A}_1
\end{equation}
where
\begin{eqnarray}
{\mathsf A}_1= \frac{L^2\tan{\theta_{_{m}}}}{4m}\equiv
A_{_{max-total}}^{^{(m>n)}}\nonumber\\
{\mathsf
A}_2=\frac{L^2\tan{\theta_{_{m}}}\tan{\theta_{_{n}}}}{4(m\tan{\theta_{_{n}}}+n\tan{\theta_{_{m}}})}
\equiv A_{_{min-total}}\end{eqnarray}

For completeness, we state the main results of this section (Eqs.
(14), (21), (22)) for the case identical polygons:
\begin{equation}
x\in\left(0,\frac{1}{2}\left(L-\sqrt{\frac{8n{\mathsf
A}-L^2\tan{\theta_{_{n}}}}{\tan{\theta_{_{n}}}}}\right)\right)\bigcup\left(\frac{1}{2}\left(L+\sqrt{\frac{8n{\mathsf
A}-L^2\tan{\theta_{_{n}}}}{\tan{\theta_{_{n}}}}}\right),L\right)
\end{equation}
subject to the constraint
\begin{equation}
2\sqrt{\frac{n{\mathsf A}}{\tan{\theta_{_{n}}}}}\leq L \leq
2\sqrt{\frac{2n{\mathsf A}}{\tan{\theta_{_{n}}}}}
\end{equation}
or equivalently
\begin{eqnarray}
\frac{L^2\tan{\theta_{_{n}}}}{8n}\leq {\mathsf A} \leq
\frac{L^2\tan{\theta_{_{n}}}}{4n}~~~~~~~~~~~~~~~~~~~~~~\nonumber\\
\textnormal{i.e.,}~~~~~~~~~~~~~~~~~~~~~~~~~~~A_{_{min-total}}^{^{(m=n)}}\leq
{\mathsf A}\leq A_{_{max-total}}^{^{(m=n)}}~~~~~~~~~~~~~~~~~~~~~~
\end{eqnarray}
{\bf\textsf{V.~~Upper bound on the area of two polygons}}\\
Now, the problem at hand is as follows:
\begin{equation}
 \frac{x^2}{4n}\tan{\theta_{_{n}}} +
\frac{(L-x)^2}{4m}\tan{\theta_{_{m}}}<{\mathsf A}
\end{equation}
Clearly, the solution to this inequality is
\begin{equation}
x\in(x_{_{min}}-\hat{x},x_{_{min}}+\hat{x})
\end{equation}
with the same constraint as before (see Eqs. 21 and 22)\\
\\
\\
{\bf\textsf{VI.~~Examples}}\\
In this sub-section we will work out examples that will illustrate
the results established in the previous sub-sections.
\\

{\bf(a.)} \texttt{Minimum area problem}

 ~~~~~~($\alpha$.)~~~\textit{hexagon} ($n=6,
~\theta_{_{n}}=\frac{\pi}{3})$ and \textit{square}
($m=4,~\theta_{_{m}}=\frac{\pi}{4})$
\begin{eqnarray}
x_{_{min}} = (2\sqrt {3}-3)L\nonumber\\
A_{_{min-total}}=\frac{(2-\sqrt{3})L^2}{8}
\end{eqnarray}
~~~~~~~~~~($\beta$.)~~\textit{circle }($n\rightarrow\infty,
~\theta_{_{n}}=\frac{\pi}{2}$ (see appendix A)) and
\textit{square}($m=4,~\theta_{_{m}}=\frac{\pi}{4})$
\begin{eqnarray}
x_{_{min}}=\frac{\pi L}{\pi+4}\nonumber\\
A_{_{min-total}}=\frac{L^2}{4(\pi+4)}
\end{eqnarray}
\\

{\bf(b.)} \texttt{Maximum area problem}

 ~~~~~~($\alpha$.)~~~\textit{hexagon}($n=6,
~\theta_{_{n}}=\frac{\pi}{3})$ and \textit{square}
($m=4,~\theta_{_{m}}=\frac{\pi}{4})$
\begin{eqnarray}
x_{_{max}}=L\nonumber\\
A_{_{max-total}}=\frac{L^2\sqrt 3}{24}
\end{eqnarray}
~~~~~~~~~~~($\beta$.)~~~\textit{pentagon} and \textit{pentagon}
($m=n=5, ~\theta_{_{m}}=\theta_{_{n}}=\frac{3\pi}{10})$
\begin{eqnarray}
x_{_{max}}=0,L\nonumber\\
A_{_{max-total}}=\frac{L^2\sqrt{5(5+2\sqrt 5)}}{100}
\end{eqnarray}
~~~~~~~~~~~~~~~~~~Note: $\tan {\frac{3\pi}{10}}=
\sqrt{\frac{5+2\sqrt 5}{5}}$
\\

 {\bf(c.)} \texttt{Lower bound problem}

  ~~~~~~($\alpha$.)~~~\textit{hexagon} ($n=6,
~\theta_{_{n}}=\frac{\pi}{3})$ and \textit{square}
($m=4,~\theta_{_{m}}=\frac{\pi}{4})$
\begin{eqnarray}
~~~~~~~~~~~x\in\left(0, \frac{3L - \sqrt{6[8(3 +2\sqrt{3}){\mathsf
A}-L^2\sqrt{3}]}}{3+ 2\sqrt{3}}\right)\bigcup \left(\frac{3L
+\sqrt{6[8(3 +2\sqrt{3}){\mathsf A}-L^2\sqrt{3}]}}{3+
2\sqrt{3}},L\right)
 \nonumber\\
4\sqrt{\mathsf A}\leq L \leq 2\sqrt{2{\mathsf
A}(2+\sqrt{3})}~~~~~~~~~~~~~~~~~~~~~~~~~~~~~~~~~~~ \nonumber\\
\frac{(2-\sqrt 3)L^2}{8}\leq {\mathsf A}\leq
\frac{L^2}{16}~~~~~~~~~~~~~~~~~~~~~~~~~~~~~~~~~~~~
\end{eqnarray}

~~~~~~($\beta$.)~~~\textit{dodecagon}($n=12,
\theta_{_{n}}=\frac{5\pi}{12})$ and \textit{triangle}
($m=3,\theta_{_{m}}=\frac{\pi}{6})$
\begin{eqnarray}
x\in\left(0, \frac{4L -2 \sqrt{12(6 +7\sqrt{3}){\mathsf
A}-L^2(3+2\sqrt{3})}}{7+ 2\sqrt{3}}\right)\bigcup \left(\frac{4L
+2 \sqrt{12(6 +7\sqrt{3}){\mathsf A}-L^2(3+2\sqrt{3})}}{7+
2\sqrt{3}},L\right)
 \nonumber\\
2\sqrt{3{\mathsf A}\sqrt{3}}\leq L \leq
2\sqrt{3{\mathsf A}(8-3\sqrt{3})}~~~~~~~~~~~~~~~~~~~~~~~~~~~~~~~~~~~ \nonumber\\
\frac{L^2(8+3\sqrt{3})}{444}\leq {\mathsf A}\leq
\frac{L^2\sqrt{3}}{36}~~~~~~~~~~~~~~~~~~~~~~~~~~~~~~~~~~~~
\end{eqnarray}
~~~~~~~~~~~~~~~~~~Note: $\tan {\frac{5\pi}{12}}= 2+\sqrt 3$
\\

 {\bf(d.)} \texttt{Upper bound problem}

  ~~~~~($\alpha$.)~~~\textit{hexagon }($n=6,
~\theta_{_{n}}=\frac{\pi}{3})$ and \textit{square}
($m=4,~\theta_{_{m}}=\frac{\pi}{4})$
\begin{eqnarray}
~~~~~~~~~~~~~~~~~~x\in\left(\frac{3L -\sqrt{6[8(3
+2\sqrt{3}){\mathsf A}-L^2\sqrt{3}]}}{3+ 2\sqrt{3}},\frac{3L
+\sqrt{6[8(3 +2\sqrt{3}){\mathsf A}-L^2\sqrt{3}]}}{3+
2\sqrt{3}}\right)
 \nonumber\\
4\sqrt{\mathsf A}\leq L \leq 2\sqrt{2{\mathsf
A}(2+\sqrt{3})}~~~~~~~~~~~~~~~~~~~~~~~~~~~~~~~~~~~ \nonumber\\
\frac{(2-\sqrt 3)L^2}{8}\leq {\mathsf A}\leq
\frac{L^2}{16}~~~~~~~~~~~~~~~~~~~~~~~~~~~~~~~~~~~~
\end{eqnarray}

~~~~~($\beta$.)~~~\textit{circle}($n\rightarrow\infty,
\theta_{_{n}}=\frac{\pi}{2})$ and \textit{octagon}
($m=8,\theta_{_{m}}=\frac{3\pi}{8})$
\begin{eqnarray}
~~~~~~~~~~~~~~~~x\in\left(\frac{8L
-\sqrt{8\pi[4(8+\pi(1+\sqrt{2})){\mathsf A}-L^2]}}{8+
\pi(1+\sqrt{2})}, \frac{8L
-\sqrt{8\pi[4(8+\pi(1+\sqrt{2})){\mathsf A}-L^2]}}{8+
\pi(1+\sqrt{2})}\right)
 \nonumber\\
2\sqrt{\pi\mathsf A}\leq L \leq 2\sqrt{{\mathsf
A}[8(\sqrt{2}-1)+\pi]}~~~~~~~~~~~~~~~~~~~~~~~~~~~~~~~~~~~ \nonumber\\
\frac{L^2}{4[8(\sqrt {2}-1)+\pi]}\leq {\mathsf A}\leq
\frac{L^2}{4\pi}~~~~~~~~~~~~~~~~~~~~~~~~~~~~~~~~~~~~
\end{eqnarray}
~~~~~~~~~~~~~~~~~~Note: $\tan {\frac{3\pi}{8}}= 1+\sqrt 2$
\\

{\bf(e.)} \texttt{An explicit numerical example}

~~~~~~Let the length of the wire be 12 units ($=L$). We would
  like to analyze the

  ~~~~~~above four problems in the case of a square ($n=4$) and a
  triangle ($m=3$).

~~~~~Firstly, to minimize the area of the polygons, the wire
should be severed

~~~~~at $x_{_{min}}=\frac{48}{4+3\sqrt 3}\approx 5.220$ units so
that the minimum total area of the

~~~~~polygons is $A_{_{min-total}}=\frac{36}{4+3\sqrt 3}\approx
3.915$ square units. This area is allocated

~~~~~between the square and triangle as
$A_{_{4}}|_{x=5.220}=\frac{144}{43+24\sqrt 3}\approx 1.703$ square

~~~~~units and $A_{_{3}}|_{x=5.220}=\frac{108\sqrt 3}{43+24\sqrt
3}\approx 2.212$ square units.

~~~~~~~Secondly, the total maximum value of the area of the
polygons is $A_{_{max-total}}=$

~~~~~$9$ square units, which means that all of the wire is used in
the perimeter of

~~~~~the square.

~~~~~For the inequality problems we first find the constraint:
$\frac{36}{4+3\sqrt 3}\leq{\mathsf A}\leq 4\sqrt 3$.

~~~~~We choose ${\mathsf A}= 5$ square units. Then,
$x_{_{min}}\pm\hat{x}=\frac{48\pm 4\sqrt{3(45-16\sqrt
3)}}{4+3\sqrt 3}\approx 8.352$,

~~~~~$2.087$ units. Thus, in the case where we want the total area
of the two

~~~~~polygons to exceed $5$ square units the wire can be cut
 either in the open

 ~~~~~interval $(0,2.087)$ or $(8.352,12)$ units. It is
 interesting to note that the

 ~~~~~corresponding areas of the square
 and triangle vary as $A_{_{4}}\in (0,0.272)\cup$

 ~~~~~$(4.360,5)$ square units and $A_{_{3}}\in (5,4.728)\cup(0.640,0)$
square units. For the

~~~~~upper bound problem, the cutting range lies in the bounded
open

~~~~~interval $(2.087,8.352)$ units.

\section{3-partitioned wire}
{\bf\textsf{I.~~The total area of three polygons}}\\
In this sub-section we consider the extrema problem when the wire
is cut at two places; the first of length, $x$, and the second
piece of length, $y$. The third segment will obviously be of
length $L-x-y\equiv z$. The total area of the $n$-, $m$- and
$p$-sided polygons is
\begin{equation}
A_{_{total}}(x,y,z) =
A_{_{n}}(x)+A_{_{m}}(y)+A_{_{p}}(z)=\frac{x^2}{4n}\tan{\theta_{_{n}}}
+\frac{y^2}{4m}\tan{\theta_{_{m}}}
+\frac{z^2}{4p}\tan{\theta_{_{p}}}
\end{equation}
In the case of identical polygons,
\begin{equation}
A_{_{total}}^{^{(m=n=p)}}=\frac{2x^2+2y^2-2L(x+y)+2xy+L^2}{4n}\tan{\theta_{_{n}}}
\end{equation}
\\
{\bf\textsf{II.~~Minimum total area of three polygons}}\\
For critical values,
\begin{eqnarray}
\frac{\partial A_{_{total}}}{\partial
x}=x\left(\frac{\tan{\theta_{_{n}}}}{2n}+\frac{\tan{\theta_{_{p}}}}{2p}\right)
+y\frac{\tan{\theta_{_{p}}}}{2p} -L\frac{\tan{\theta_{_{p}}}}{2p}
=0\nonumber\\
\frac{\partial A_{_{total}}}{\partial
y}=x\frac{\tan{\theta_{_{p}}}}{2p}+y\left(\frac{\tan{\theta_{_{m}}}}{2m}+\frac{\tan{\theta_{_{p}}}}{2p}\right)
 -L\frac{\tan{\theta_{_{p}}}}{2p}
=0
\end{eqnarray}
Solving these two simultaneous linear equations for $x$ and $y$
results in the minimum value for $A_{_{total}}$:
\begin{equation}
x_{_{min}}=\frac{nL\tan{\theta_{_{m}}}\tan{\theta_{_{p}}}}
{n\tan{\theta_{_{m}}}\tan{\theta_{_{p}}}+m\tan{\theta_{_{n}}}\tan{\theta_{_{p}}}
+p\tan{\theta_{_{m}}}\tan{\theta_{_{n}}}}
\end{equation}
\begin{equation}
y_{_{min}}=\frac{mL\tan{\theta_{_{n}}}\tan{\theta_{_{p}}}}
{n\tan{\theta_{_{m}}}\tan{\theta_{_{p}}}+m\tan{\theta_{_{n}}}\tan{\theta_{_{p}}}
+p\tan{\theta_{_{m}}}\tan{\theta_{_{n}}}}
\end{equation}
because
\begin{eqnarray}
\frac{\partial^2A_{_{total}}}{\partial x^2}\frac{\partial^2
A_{_{total}}}{\partial y^2}-\left( \frac{\partial^2
A_{_{total}}}{\partial x\partial
y}\right)^2=\frac{n\tan{\theta_{_{m}}}\tan{\theta_{_{p}}}+m\tan{\theta_{_{n}}}\tan{\theta_{_{p}}}
+p\tan{\theta_{_{m}}}\tan{\theta_{_{n}}}}{4nmp}>0 \nonumber\\
\textnormal{and}~~~~~
 \frac{\partial^2A_{_{total}}}{\partial
x^2}+\frac{\partial^2 A_{_{total}}}{\partial
y^2}=\frac{\tan{\theta_{_{n}}}}{2n}+\frac{\tan{\theta_{_{m}}}}{2m}+\frac{\tan{\theta_{_{p}}}}{p}>0~~~~~~~~~~~~~~~~~~~~~~~~~~~~~~~
\end{eqnarray}
Further,
\begin{equation}
z_{_{min}}=\frac{pL\tan{\theta_{_{n}}}\tan{\theta_{_{m}}}}
{n\tan{\theta_{_{m}}}\tan{\theta_{_{p}}}+m\tan{\theta_{_{n}}}\tan{\theta_{_{p}}}
+p\tan{\theta_{_{m}}}\tan{\theta_{_{n}}}}
\end{equation}
 Therefore,
\begin{eqnarray}
A_{_{n}}(x_{_{min}})=\frac{L^2n\tan^2{\theta_{_{m}}}\tan^2{\theta_{_{p}}}\tan{\theta_{_{n}}}}
{4\left(n\tan{\theta_{_{m}}}\tan{\theta_{_{p}}}+m\tan{\theta_{_{n}}}\tan{\theta_{_{p}}}
+p\tan{\theta_{_{m}}}\tan{\theta_{_{n}}}\right)^2}\nonumber\\
A_{_{m}}(y_{_{min}})=\frac{L^2m\tan^2{\theta_{_{n}}}\tan^2{\theta_{_{p}}}\tan{\theta_{_{m}}}}
{4\left(n\tan{\theta_{_{m}}}\tan{\theta_{_{p}}}+m\tan{\theta_{_{n}}}\tan{\theta_{_{p}}}
+p\tan{\theta_{_{m}}}\tan{\theta_{_{n}}}\right)^2}\nonumber\\
A_{_{p}}(z_{_{min}})=\frac{L^2p\tan^2{\theta_{_{n}}}\tan^2{\theta_{_{m}}}\tan{\theta_{_{p}}}}
{4\left(n\tan{\theta_{_{m}}}\tan{\theta_{_{p}}}+m\tan{\theta_{_{n}}}\tan{\theta_{_{p}}}
+p\tan{\theta_{_{m}}}\tan{\theta_{_{n}}}\right)^2}\nonumber\\
A_{_{min-total}} \equiv
A_{_{total}}(x_{_{min}},y_{_{min}},z_{_{min}})
=\frac{L^2\tan{\theta_{_{n}}}\tan{\theta_{_{m}}}\tan{\theta_{_{p}}}}
{4\left(n\tan{\theta_{_{m}}}\tan{\theta_{_{p}}}+m\tan{\theta_{_{n}}}\tan{\theta_{_{p}}}
+p\tan{\theta_{_{m}}}\tan{\theta_{_{n}}}\right)}
\end{eqnarray}

In the case of identical polygons
\begin{eqnarray}
x^{^{(m=n=p)}}_{_{min}}=y^{^{(m=n=p)}}_{_{min}}=z^{^{(m=n=p)}}_{_{min}}=\frac{L}{3}\nonumber\\
A_{_{n}}(x_{_{min}})=A_{_{m}}(y_{_{min}})=A_{_{p}}(z_{_{min}})=
\frac{L^2}{36n}\tan{\theta_{_{n}}}\nonumber\\
A_{_{min-total}}^{^{(m=n=p)}}= \frac{L^2}{12n}\tan{\theta_{_{n}}}
\end{eqnarray}
\\
{\bf\textsf{III.~~Maximum total area of three polygons}}\\
To compute the maximum total area, we must consider the boundary
values at $x=0$ or $y=0$ or $z=0$. We begin by finding the
critical values $y_{_{cr}}$ and $z_{_{cr}}$ of
$A_{_{total}}(0,y,z)$. This total area, when evaluated at these
critical values, becomes a candidate for the maximum total area.
The results are
\begin{eqnarray}
A_{_{total}}(0,y_{_{cr}},z_{_{cr}})=
\frac{L^2\tan{\theta_{_{p}}}\tan{\theta_{_{m}}}}{4(p\tan{\theta_{_{m}}}+
m\tan{\theta_{_{p}}})}\nonumber\\
y_{_{cr}} = \frac{mL\tan{\theta_{_{p}}}}{m\tan{\theta_{_{p}}}+
p\tan{\theta_{_{m}}}}\nonumber\\
z_{_{cr}} = \frac{pL\tan{\theta_{_{m}}}}{m\tan{\theta_{_{p}}}+
p\tan{\theta_{_{m}}}}
\end{eqnarray}
Similarly,
\begin{eqnarray}
A_{_{total}}(x_{_{cr}},0,z_{_{cr}}')=
\frac{L^2\tan{\theta_{_{n}}}\tan{\theta_{_{p}}}}{4(n\tan{\theta_{_{p}}}+
p\tan{\theta_{_{n}}})}\nonumber\\
x_{_{cr}} = \frac{nL\tan{\theta_{_{p}}}}{n\tan{\theta_{_{p}}}+
p\tan{\theta_{_{n}}}}\nonumber\\
z_{_{cr}}' = \frac{pL\tan{\theta_{_{n}}}}{n\tan{\theta_{_{p}}}+
p\tan{\theta_{_{n}}}}
\end{eqnarray}
and
\begin{eqnarray}
A_{_{total}}(x_{_{cr}}',y_{_{cr}}',0)=
\frac{L^2\tan{\theta_{_{n}}}\tan{\theta_{_{m}}}}{4(n\tan{\theta_{_{m}}}+
m\tan{\theta_{_{n}}})}\nonumber\\
x_{_{cr}}' = \frac{nL\tan{\theta_{_{m}}}}{n\tan{\theta_{_{m}}}+
m\tan{\theta_{_{n}}}}\nonumber\\
y_{_{cr}}' = \frac{mL\tan{\theta_{_{n}}}}{n\tan{\theta_{_{m}}}+
m\tan{\theta_{_{n}}}}
\end{eqnarray}
Finally, the maximum total area is given by
\begin{equation}
A_{_{max-total}}=\max\left\{A_{_{total}}(0,y_{_{cr}},z_{_{cr}}),~
A_{_{total}}(x_{_{cr}},0,z_{_{cr}}'),~
A_{_{total}}(x_{_{cr}}',y_{_{cr}}',0)\right\}
\end{equation}
\\
{\bf\textsf{IV.~~Upper and lower bound on the area of three
polygons}}
\\
We consider a simple inequality problem in which the first two
polygons have the same perimeter $x$. As before, we look for
values of $x$ for which the combined area of these polygons is
greater than or less than a specified area, ${\mathsf A}$.  Thus,
for the lower bound problem, we can write
\begin{equation}
 \frac{x^2}{4n}\tan{\theta_{_{n}}} +\frac{x^2}{4m}\tan{\theta_{_{m}}}
+\frac{(L-2x)^2}{4p}\tan{\theta_{_{p}}}>{\mathsf A}
\end{equation}
The solution to this inequality is given by
\begin{equation}
x\in(0,x_-) \cup (x_+,L)
\end{equation}
where
\begin{eqnarray}
x_{\pm}=\frac{1}{{mp\tan{\theta_{_{n}}}+np\tan{\theta_{_{m}}}+4mn\tan{\theta_{_{p}}}}}\times
~~~~~~~~~~~~~~~~~~~~~~~~~~~~~~\nonumber\\
\{2Lmn\tan{\theta_{_{p}}}
 \pm (mnp[4{\mathsf A}(mp\tan{\theta_{_{n}}}+np\tan{\theta_{_{m}}}+4mn\tan{\theta_{_{p}}})\nonumber\\
-L^2\tan{\theta_{_{p}}}(m\tan{\theta_{_{n}}}+n\tan{\theta_{_{m}}})])^{\frac{1}{2}}\}
\end{eqnarray}
Obviously, the solution to the upper bound problem is
\begin{equation}
x\in(x_-, x_+)
\end{equation}
The above solutions are valid provided that
\begin{equation}
\frac{L^2\tan{\theta_{_{p}}}(m\tan{\theta_{_{n}}}+n\tan{\theta_{_{m}}})}{4(mp\tan{\theta_{_{n}}}
+np\tan{\theta_{_{m}}}+4mn\tan{\theta_{_{p}}})}\leq {\mathsf
A}\leq \frac{L^2\tan{\theta_{_{p}}}}{4p}
\end{equation}
\\
{\bf\textsf{V.~~Example}}\\
A wire 10 units ($=L$) long is cut into three pieces. The first
piece is bent to form a square ($n=4$), the second to form a
triangle ($m=3$), and the third to form a circle
($p\rightarrow\infty$). We are interested in knowing where the
cuts should be made so that the sum of the three areas is a
minimum and also the manner of cutting that
will produce a maximum total area.\\
For the minima problem, the length of the first piece is
$x_{_{min}}=\frac{40}{4+3\sqrt 3+\pi}\approx 3.242$ units, the
second piece is of length $y_{_{min}}=\frac{30\sqrt 3}{4+3\sqrt
3+\pi}\approx 4.212$ units, and the third piece automatically
comes out to be $z_{_{min}}=\frac{10\pi}{4+3\sqrt 3+\pi}\approx
2.546$ units. This gives a total minimum area of
$A_{_{min-total}}=\frac{25}{4+3\sqrt 3+\pi}\approx 2.026$ square
units with the assignment
$A_{_{4}}|_{x=3.242}=\frac{100}{(4+3\sqrt 3+\pi)^2}\approx 0.657$
square units, $A_{_{3}}|_{y=4.212}=\frac{75\sqrt 3}{(4+3\sqrt
3+\pi)^2}\approx 0.853$ square units and
$A_{_{\infty}}|_{z=2.546}=\frac{25\pi}{(4+3\sqrt 3+\pi)^2}\approx
0.516$ square units.\\
The maximum total area of the polygons is
$A_{_{max-total}}=\frac{25\pi}{4+\pi}\approx 10.998$ square units
and occurs at the boundary $y=0$ with critical values
$x_{_{c}}=\frac{40}{4+\pi}\approx 5.601$ units and
$z_{_{c}}=\frac{10\pi}{4+\pi}\approx 4.399$ units.\\
In the case of the inequality problem, the area ${\mathsf A}$ must
be chosen so that \\ $\frac{25(4+3\sqrt 3)}{3\pi\sqrt
3+4\pi+48\sqrt 3}\leq {\mathsf A}\leq \frac{25}{\pi}$ square
units. We choose ${\mathsf A}$=5 square units. This gives
$x_{\pm}=\frac{240\sqrt 3\pm \sqrt{240\pi(9\pi+4\pi\sqrt
3+99-20\sqrt 3)}}{3\pi\sqrt 3+4\pi+48\sqrt 3}\approx 6.332, 1.089$
units. Thus, in the case of the upper bound problem $x\in
(1.089,6.332)$ units while for the lower bound problem $x\in
(0,1.089)\cup (6.332,10)$.

\section{(${\mathsf k}$\bf {+1})-partitioned wire}
{\bf\textsf{I.~~The total area of the (k+1)-polygons}}\\
In this sub-section we consider the extrema problem where the wire
is cut at $\mathsf {k}$ different places with lengths, $x_{_1}$,
$x_{_2}$,...., $x_{_{\mathsf k}}$. The last segment being of
length, $x_{_{\mathsf k+1}}=L-x_{_1}-x_{_2}-.....-x_{_{\mathsf
k}}$. Then, the total area of the $n_{_1}$-, $n_{_2}$-,.....,
$n_{_{\mathsf k+1}}$-sided polygons is

\begin{eqnarray}
A_{_{total}}(x_{_1},x_{_2},...,x_{_{\mathsf k+1}})
=\sum_{i=1}^{{\mathsf k}+1}A_{_{n_{_i}}}(x_{_i})=
\frac{x_{_1}^2}{4n_{_1}}\tan{\theta_{_{n_{_1}}}}
+\frac{x_{_2}}{4n_{_2}}\tan{\theta_{_{n_{_2}}}}
+.....+\frac{x_{_{\mathsf k+1}}^2}{4n_{_{\mathsf
k+1}}}\tan{\theta_{_{n_{_{\mathsf
k+1}}}}}\nonumber\\
=\sum_{i=1}^{{\mathsf
k}+1}\frac{x_{_i}^2}{4n_{_i}}\tan{\theta_{_{n_{_i}}}}~~~~~~~~~~~~~~~~~~~~~~~~~~~~~~~~~~~~~~~~~~~~~
\end{eqnarray}
In the case of identical polygons,
\begin{equation}
A_{_{total}}^{^{(n_{_1}=n_{_2}=...=n_{_{\mathsf k+1}}\equiv n)}}
=\frac{\tan{\theta_{_{n}}}}{4n}\left[L^2+2\sum_{i=1,i\neq
j}^{\mathsf k}\left(x_{_i}^2+x_{_i}x_{_j}-Lx_{_i}\right) \right]
\end{equation}
\\
{\bf\textsf{II.~~Minimum total area of (k+1)-polygons}}\\
For critical values
\begin{equation}
\frac{\partial A_{_{total}}}{\partial x_{_1}}=\frac{\partial
A_{_{total}}}{\partial x_{_2}}=.....=\frac{\partial
A_{_{total}}}{\partial x_{_{\mathsf k}}}=0
\end{equation}
with
\begin{equation}
\frac{\partial A_{_{total}}}{\partial x_{_i}}= x_{_i}\left(
\frac{\tan{\theta_{_{n_{_i}}}}}{2n_{_i}}+\frac{\tan{\theta_{_{n_{_{\mathsf
k+1}}}}}}{{2n_{_{\mathsf k+1}}}}
\right)+\frac{\tan{\theta_{_{n_{_{\mathsf k+1}}}}}}{{2n_{_{\mathsf
k+1}}}} \left(-L
 +\sum_{j=1,j\ne i}^{\mathsf k}x_{_j}\right);~~~(i=1~
{\textnormal{to}}~{\mathsf k})
\end{equation}

 Solving these ${\mathsf k}$ simultaneous linear equations for $x_{_1}$,
$x_{_2}$,...., $x_{_{\mathsf k}}$ results in the minimum value for
$A_{_{total}}$:
\begin{equation}
x_{_{i_{_{min}}}}=L\left(\frac{n_{_i}}{\tan{\theta_{_{n_{_i}}}}}\right)\left(\sum_{j=1}^{\mathsf
k+1}\frac{n_{_j}}{\tan{\theta_{_{n_{_j}}}}}\right)^{-1}=L\left(\frac{n_{_i}\prod_{j=1,j\ne
i}^{\mathsf k+1}\tan{\theta_{_{n_{_j}}}}} {\sum_{l=1}^{\mathsf
k+1}n_{_l}\prod_{j=1,j\ne l}^{\mathsf
k+1}\tan{\theta_{_{n_{_j}}}}}\right);~~~(i=1~
{\textnormal{to}}~{\mathsf k})
\end{equation}

because
\begin{eqnarray}
\det{\left(\matrix{\frac{\partial^2A_{_{total}}}{\partial x_{_1}
^2} & \frac{\partial^2 A_{_{total}}}{\partial x_{_1}\partial
x_{_2}} &.&.&.&.&\frac{\partial^2 A_{_{total}}}{\partial
x_{_1}\partial x_{_{\mathsf k}}}\cr
 \frac{\partial^2 A_{_{total}}}{\partial
x_{_2}\partial x_{_1}} & \frac{\partial^2A_{_{total}}}{\partial
x_{_2} ^2}&.&.&.&.& \frac{\partial^2 A_{_{total}}}{\partial
x_{_2}\partial x_{_{\mathsf k}}} \cr .&&&&&&. \cr .&&&&&&. \cr
.&&&&&&.\cr .&&&&&&. \cr \frac{\partial^2 A_{_{total}}}{\partial
x_{_{\mathsf k}}\partial x_{_1}}&.&.&.&.&.&\frac{\partial^2
A_{_{total}}}{\partial x_{_{\mathsf k}}^2}} \right)}>0\nonumber\\
\textnormal{and}~~~~~~~~~~~~~~~~~~~~~~~~~~~~~~~~~~~~~~~~~~~~~~~~~~~~\sum_{i=1}^{\mathsf
k}\frac{\partial^2A_{_{total}}}{\partial x_{_i}
^2}>0~~~~~~~~~~~~~~
\end{eqnarray}
Therefore,
\begin{eqnarray}
A_{_{n_{_i}}}(x_{_{i_{_{min}}}})=\frac{L^2}{4}\left(\frac{n_{_i}\tan{\theta_{_{n_{_i}}}}\prod_{j=1,j\ne
i}^{\mathsf k}\tan^2{\theta_{_{n_{_j}}}}}
{\left[\sum_{l=1}^{\mathsf k+1}n_{_l}\prod_{j=1,j\ne l}^{\mathsf
k+1}\tan{\theta_{_{n_{_j}}}}\right]^2}\right);~~~(i=1~
{\textnormal{to}}~{\mathsf k})\nonumber\\
A_{_{min-total}}=\frac{L^2}{4}\left(\sum_{j=1}^{\mathsf
k+1}\frac{n_{_j}}{\tan{\theta_{_{n_{_j}}}}}\right)^{-1}=\frac{L^2}{4}\left(\frac{\prod_{j=1}^{\mathsf
k+1}\tan{\theta_{_{n_{_j}}}}} {\sum_{l=1}^{\mathsf
k+1}n_{_l}\prod_{j=1,j\ne l}^{\mathsf
k+1}\tan{\theta_{_{n_{_j}}}}}\right)
\end{eqnarray}

In the case of identical polygons
\begin{eqnarray}
x^{^{(n_{_1}=n_{_2}=...=n_{_{\mathsf
k+1}})}}_{_{i_{_{min}}}}=\frac{L}{{\mathsf k}+1};~~~(i=1~
{\textnormal{to}}~{\mathsf k})\nonumber\\
A_{_{n_{_i}}}^{^{(n_{_1}=n_{_2}=...=n_{_{\mathsf k+1}}\equiv
n)}}|_{{_{x=x_{_{i_{_{min}}}}}}}= \frac{L^2}{4({\mathsf
k}+1)^2n}\tan{\theta_{_{n}}};~~~(i=1~
{\textnormal{to}}~{\mathsf k})\nonumber\\
A_{_{min-total}}^{^{(n_{_1}=n_{_2}=...=n_{_{\mathsf k+1}}\equiv
n)}}= \frac{L^2}{4({\mathsf k}+1)n}\tan{\theta_{_{n}}}~~~~~~~~~~
\end{eqnarray}
\\
{\bf\textsf{III.~~Maximum total area of (k+1)-polygons}}\\
Here we need to consider the total area at the boundary values,
$x_{_b}=0;~b=1$ to ${\mathsf k+1}$:
\begin{equation}
A_{_{total}}(x_{_1},x_{_2},..,x_{_b}=0,..,x_{_{\mathsf k+1}})
=\sum_{i=1,i\neq b}^{{\mathsf
k}+1}\frac{x_{_i}^2}{4n_{_i}}\tan{\theta_{_{n_{_i}}}};~~~(b=1~
{\textnormal{to}}~{\mathsf k}+1)
\end{equation}

The critical value, $x_{_{i_{_{cr}}}}$, of this $A_{_{total}}$ is
given by
\begin{equation}
x_{_{i_{_{cr}}}}=L\left(\frac{n_{_i}\prod_{j=1,j\ne i,b}^{\mathsf
k+1}\tan{\theta_{_{n_{_j}}}}} {\sum_{l=1,l\ne b}^{\mathsf
k+1}n_{_l}\prod_{j=1,j\ne l,b}^{\mathsf
k+1}\tan{\theta_{_{n_{_j}}}}}\right);~~~(i=1~
{\textnormal{to}}~{\mathsf k}~{\textnormal{except}}~b)
\end{equation}

The above $A_{_{total}}$, when evaluated at these critical values,
becomes a candidate for the maximum total area:
\begin{equation}
A_{_{total}}(x_{_{1_{_{cr}}}},x_{_{2_{_{cr}}}},..,x_{_b}=0,..,x_{_{{\mathsf
({k}+1)}_{_{cr}}}})=\frac{L^2}{4}\left(\frac{\prod_{j=1,j\neq
b}^{\mathsf k+1}\tan{\theta_{_{n_{_j}}}}} {\sum_{l=1,l\neq
b}^{\mathsf k+1}n_{_l}\prod_{j=1,j\ne l,b}^{\mathsf
k+1}\tan{\theta_{_{n_{_j}}}}}\right);~~~(b=1~
{\textnormal{to}}~{\mathsf k}+1)
\end{equation}
and
\begin{eqnarray}
 A_{_{max-total}}=\max_{b=1}^{{\mathsf
k}+1}\left\{A_{_{total}}(x_{_{1_{_{cr}}}},x_{_{2_{_{cr}}}},..,x_{_b}=0,..,x_{_{{\mathsf
({k}+1)}_{_{cr}}}})\right\}
\end{eqnarray}
\\
{\bf\textsf{IV.~~Upper and lower bound on the area of (k+1)-polygons}}\\
As in the 2-cut case, we consider the simple case in which the
first ${\mathsf k}$ polygons have the same perimeter, $x$.
Therefore, the last one will have a perimeter of $L-{\mathsf k}x$.
Hence, for the lower bound problem
\begin{equation}
\frac{x^2}{4n_{_1}}\tan{\theta_{_{n_{_1}}}}
+\frac{x^2}{4n_{_2}}\tan{\theta_{_{n_{_2}}}}
+.....+\frac{x^2}{4n_{_{\mathsf k}}}\tan{\theta_{_{n_{_{\mathsf
k}}}}}+\frac{(L-{\mathsf k}x)^2}{{4n_{_{\mathsf
k+1}}}}\tan{\theta_{_{n_{_{\mathsf k+1}}}}}>{\mathsf A}
\end{equation}
This can be simplified to
\begin{eqnarray}
x^2\left[\sum_{l=1}^{\mathsf
k}\tan{\theta_{_{n_{_l}}}}\prod_{j=1,j\ne l}^{\mathsf
k+1}n_{_j}+{\mathsf k}^2 \tan{\theta_{_{n_{_{\mathsf
k+1}}}}}\prod_{j=1}^{\mathsf k}n_{_j}\right]-2L{\mathsf
k}x\tan{\theta_{_{n_{_{\mathsf k+1}}}}}\prod_{j=1}^{\mathsf
k}n_{_j}\nonumber\\
+L^2\tan{\theta_{_{n_{_{\mathsf k+1}}}}}\prod_{j=1}^{\mathsf
k}n_{_j}-4{\mathsf A}\prod_{j=1}^{\mathsf k+1}n_{_j}>0 \nonumber\\
\end{eqnarray}
The solution is
\begin{equation}
x\in (0,x_-)\cup (x_+,L)
\end{equation}
where
\begin{eqnarray}
x_{\pm}=\frac{L{\mathsf k}\tan{\theta_{_{n_{_{\mathsf
k+1}}}}}\prod_{j=1}^{\mathsf k}n_{_j}\pm
\sqrt{\prod_{j=1}^{\mathsf k+1}n_{_j}\sum_{l=1}^{{\mathsf
k+1}}(4{\mathsf A}\alpha_l-L^2\tan{\theta_{_{n_{_{\mathsf
k+1}}}}}\beta_l) \tan{\theta_{_{n_{_l}}}}\prod_{i=1,i\ne
l}^{\mathsf k+1}n_{_i}}} {\sum_{l=1}^{\mathsf
k+1}\alpha_l\tan{\theta_{_{n_{_l}}}}\prod_{i=1,i\ne l}^{\mathsf
k+1}n_{_i}}\nonumber\\
\end{eqnarray}
and
\begin{eqnarray}
\alpha_{\mathsf k+1}={\mathsf k}^2,~~~~~\beta_{\mathsf
k+1}=0\nonumber\\
 \alpha_i=\beta_i=1,~~~~\forall i\ne {\mathsf k}+1
\end{eqnarray}
Finally, the above solutions are valid if
\begin{equation}
\frac{L^2\tan{\theta_{_{n_{_{\mathsf k+1}}}}}\sum_{l=1}^{\mathsf
k}\tan{\theta_{_{n_{_l}}}}\prod_{j=1,j\ne l}^{\mathsf k}n_{_j}}
{4\sum_{l=1}^{\mathsf
k+1}\alpha_l\tan{\theta_{_{n_{_l}}}}\prod_{j=1,j\ne l}^{\mathsf
k+1}n_{_j}}\leq {\mathsf A}\leq \frac{L^2
\tan{\theta_{_{n_{_{\mathsf k+1}}}}}}{4n_{_{\mathsf k+1}}}
\end{equation}
The upper bound solution is
\begin{equation}
x\in(x_-, x_+)
\end{equation}
\\
{\bf\textsf{V.~~Example}}\\
We wrap-up this section with one final concrete example. To that
end, consider a wire of length 20 units. Further, assume that the
wire is cut into six segments, with the first one bent into a
triangle ($n_{_1}=3$), the second into a square ($n_{_2}=4$), the
third into a hexagon ($n_{_3}=6$), the fourth into an octagon
($n_{_4}=8$), the fifth into a dodecagon ($n_{_5}=12$), and the
sixth into a circle ($n_{_6}\rightarrow\infty$). As before, we
consider the mannerism in which the wire should be cut so that the
total area enclosed by the above six polygons is a minimum. Here
we also deal with the
corresponding maxima problem.\\
For the minima problem, the wire should be severed at
$x_{_{1_{_{min}}}}=4.654$ units, $x_{_{2_{_{min}}}}=3.582$ units,
$x_{_{3_{_{min}}}}=3.103$ units, $x_{_{4_{_{min}}}}=2.968$ units,
$x_{_{5_{_{min}}}}=2.880$ units, and $x_{_{6_{_{min}}}}=2.814$
units giving a $A_{_{min-total}}=4.478$ square units with the
allocation $A_{_{3}}|_{x_{_1}=4.654}=1.042$ square units,
$A_{_{4}}|_{x_{_2}=3.852}=0.802$ square units,
$A_{_{6}}|_{x_{_3}=3.103}=0.695$ square units,
$A_{_{8}}|_{x_{_4}=2.968}=0.665$ square units,
$A_{_{12}}|_{x_{_5}=2.880}=0.645$ square units,
$A_{_{\infty}}|_{x_{_6}=2.814}=0.630$ square units.\\
The maximum total area of the polygons is calculated to be
$A_{_{max-total}}=5.836$ square units and occurs at the boundary
point $x_{_1}=0$ with critical values, $x_{_{2_{_{cr}}}}=4.669$
units, $x_{_{3_{_{cr}}}}=4.043$ units, $x_{_{4_{_{cr}}}}=3.868$
units, $x_{_{5_{_{cr}}}}=3.753$ units, $x_{_{6_{_{cr}}}}=3.667$
units.\\
For the inequality problem, we must choose ${\mathsf A}$ between
4.599 and 31.831 square units. Take ${\mathsf A}= 23$ square
units. Then, $x_{\pm}\approx 6.235,~0.609$ units. Thus, for the
upper bound problem $x\in (0.609, 6.235)$ units, while for lower
bound problem $x\in(0,0.609)\cup(6.235,20)$ units.

\section{Extrema problem of the distribution of a fixed number of polygonal sides among an arbitrary number of
wires} In this section we first consider two wires of specified
lengths which are bent to form two regular polygons. The total
number of sides of the polygons, $I$, is fixed. The question we
answer is how the total number of polygonal sides, $I$, should be
divided between the two wires so that the total area of the two
polygons is maximized.

To solve the problem at hand, we consider two wires of lengths
$L_n$ and $L_m$. These are bent into $n$- and $m$-sided regular
polygons with $n+m(=I)$ fixed. The area of the two polygons are
given by
$\frac{L_n^2}{4n}\tan\left(\frac{\pi}{2}-\frac{\pi}{n}\right)$ and
$\frac{L_m^2}{4(I-n)}\tan\left(\frac{\pi}{2}-\frac{\pi}{I-n}\right)$.
Thus,
\begin{equation}
A_{_{total}}=\frac{L_n^2}{4n}\tan\left(\frac{\pi}{2}-\frac{\pi}{n}\right)+
\frac{L_m^2}{4(I-n)}\tan\left(\frac{\pi}{2}-\frac{\pi}{I-n}\right)
\end{equation}
Differentiating with respect to the number of polygonal sides,
$n$, we get,
\begin{eqnarray}
\frac{dA_{_{total}}}{dn}=\frac{L_n^2}{4n^2}\left[\frac{\pi}{n}+\frac{\pi}{n}
{\textnormal{cot}}^2\left(\frac{\pi}{n}\right)-{\textnormal{cot}}\left(\frac{\pi}{n}\right)\right]
~~~~~~~~~~~~~~~~~~~~~~~~~~~~~~~~~~~~~\nonumber\\
+\frac{L_m^2}{4(I-n)^2}\left[-\frac{\pi}{I-n}-\frac{\pi}{I-n}
{\textnormal{cot}}^2\left(\frac{\pi}{I-n}\right)+{\textnormal{cot}}\left(\frac{\pi}{I-n}\right)\right]=0
\end{eqnarray}
This gives a nonlinear equation in $n$:
\begin{equation}
(\alpha L_n)^2\left[\alpha
{\textnormal{cot}}^2\alpha-{\textnormal{cot}}\alpha +\alpha\right]
=(\beta L_m)^2\left[\beta
{\textnormal{cot}}^2\beta-{\textnormal{cot}}\beta +\beta\right]
\end{equation}
where
\begin{eqnarray}
\alpha=\frac{\pi}{n}\nonumber\\
\beta=\frac{\alpha \pi}{\alpha I-n}=\frac{\pi}{I-n}
\end{eqnarray}

Finally, we generalize the above analysis for the case when there
are $\mathrm{k}$ wires with lengths $L_1$, $L_2$,....,
$L_{\mathrm{k}}$. As before, we bent each of these wires into
regular $n_{_1}$-, $n_{_2}$-,...., $n_{_{\mathrm{k}}}$-gons with
$n_{_1}+n_{_2}+....+n_{_{\mathrm{k}}}=I$. We look for values of
$n_{_1}$, $n_{_2}$,...., $n_{_{\mathrm{k}}}$ which would maximize
the total area of the polygons. In this case, the total area of
the polygons is given by
\begin{equation}
A_{_{total}}=\frac{1}{4}\sum_{i=1}^{\mathrm
{k}}\frac{L_i^2}{n_{_i}}\tan{\left(\frac{\pi}{2}-\frac{\pi}{n_{_i}}\right)}
\end{equation}
To solve the extrema problem we set
\begin{equation}
\frac{\partial A_{_{total}}}{\partial n_{_1}}=\frac{\partial
A_{_{total}}}{\partial n_{_2}}=.....=\frac{\partial
A_{_{total}}}{\partial n_{_{\mathrm{k}-1}}}=0
\end{equation}
obtaining
\begin{eqnarray}
(\alpha _1L_1)^2\left[\alpha_1
{\textnormal{cot}}^2\alpha_1-{\textnormal{cot}}\alpha_1
+\alpha_1\right]~~~~~~~~~~~~~~~~~~\nonumber\\
=(\alpha_2 L_2)^2\left[\alpha_2
{\textnormal{cot}}^2\alpha_2-{\textnormal{cot}}\alpha_2
+\alpha_2\right]~~~~~~~~~~~~~~~~~\nonumber\\
=........~~~~~~~~~~~~~~~~~~~~~~~~~~~~~~~~~~~~~~~~~~~~~~~~~~~~~\nonumber\\
=(\alpha_{\mathrm{k}-1}
L_{\mathrm{k}-1})^2\left[\alpha_{\mathrm{k}-1}
{\textnormal{cot}}^2\alpha_{\mathrm{k}-1}-{\textnormal{cot}}\alpha_{\mathrm{k}-1}
+\alpha_{\mathrm{k}-1}\right]\nonumber\\
=(\gamma L_{\mathrm{k}})^2\left[\gamma
{\textnormal{cot}}^2\gamma-{\textnormal{cot}}\gamma
+\gamma\right]~~~~~~~~~~~~~~~~~~~~~~~~
\end{eqnarray}
where
\begin{equation}
\alpha_i=\frac{\pi}{n_{_i}}~~{\textnormal{and}}~~\gamma=\frac{\pi}{I-\sum_{j=1}^{\mathrm{k}-1}n_{_j}}
\end{equation}

Appendix B contains a C program that singles out the $\{n_{_1},
n_{_2},...,n_{_k}\}$ set for which the total area (Eq.(78))
attains its maximum. Rather than solving Eqs.(80) the program
scans the entire subspace of permissible partitions ($n_{_1},
n_{_2},...,n_{_k}$;$~n_{_1}+n_{_2}+...+n_{_k}=I$) while it keeps
track of the largest value of $A_{total}$. Statements printing out
each set and the corresponding $A_{total}$ are included in the
code, but have been commented out for the sake of execution speed.

\section{Conclusions}
We began this paper with a single cut of a piece of wire and
computed the lengths of the segments formed so that the total
combined area of the arbitrary sided polygons bent from these two
segments is maximized and minimized.  After solving this extrema
problem, we went on to find the cutting range so that the combined
area of these polygons is greater than or less than a
user-specified area.  Here we also mention the constraint on the
selection of this user-specified area. For the reader we have
worked out examples in exact form, leaving the answer in terms of
the length of the wire.  These examples consist of all kinds of
combinations of polygons, however, for each of the above four
cases, a hexagon and a square is examined to get a feel of how
this combination carries itself out for each of the four cases.

We have extended the results of the above four cases to three
partitioned wires and finally to any number of partitioned wires.
However, for the inequality problem with more than one cut, we
considered the simplified case where all of the polygons are
assigned the same variable perimeter except the last one which is
automatically determined from the total length of the wire. To
illustrate to the reader the use of these results, we considered a
5-cut case (6-partitioned) of a wire of length 20 units . The
segments are bent into a triangle, square, hexagon, octagon,
dodecagon, and a circle. We found the following assignment for the
perimeter of the polygons: 4.654, 3.582, 3.103, 2.968, 2.880, and
2.814 units respectively, so that the combined area of these
polygons is minimized with this area given by 4.478 square units.
For the total area of the polygons to be maximized, we found the
following allocation for the perimeter of the polygons: 4.669 for
the square, 4.043 for the hexagon, 3.868 for the octagon, 3.753
for the dodecagon, and 3.667 for the circle. Note that for this
case none of the wire is used for the triangle. For the inequality
problem, we chose the user-specified area to be 23 square units
from a possible range of 4.599 to 31.831 square units. So, for the
upper bound problem, this result gives a cutting range between
0.609 and 6.235 units for the first five polygons and for the
lower bound problem a cutting range between 0 and 0.609 or 6.235
and 20 units.

In section 5 we explored the case where the number of sides of the
polygons formed is fixed and we calculated how this fixed number
of sides can be allocated firstly, among two wires, and then an
arbitrary number of wires so that the total area of the polygons
is maximized.

\section{Acknowledgements}
We wish to thank Health Careers Academy for allowing and
encouraging us to carry on many fruitful discussions surrounding
this paper. Two of us (RMS \& AC) wishes to also thank the
Department of Physics at Northeastern University for providing an
environment conducive to research.
\\

\section{Appendix A}
In this appendix we will show that the maximum area, ${\mathcal
A}_{_{n}}$, of an $n$-sided polygon (with perimeter, $P$) is
attained in the limiting case of a circle: $n\rightarrow\infty,
\theta_{_{n}}=\pi/2$.
\begin{eqnarray}
{\mathcal
A}_{_{n}}=\frac{P^2}{4n}\tan{\theta_{_{n}}};~~~\theta_{_{n}} =
(\frac{1}{2}-
\frac{1}{n})\pi \nonumber\\
\frac{d{\mathcal
A}_{_{n}}}{d\theta_{_{n}}}=\left[\frac{P(\pi-2\theta_{_{n}})}{4\pi}
    \right]^2\left[(\frac{\pi}{2}-\theta_{_{n}})\frac{1}{\cos^2\theta_{_{n}}}-\frac{\sin\theta_{_{n}}}{\cos\theta_{_{n}}}
    \right]=0
\end{eqnarray}
 so that
 \begin{eqnarray}
\textnormal{either}~~~~~~~~~~~~~~~~~~~~~~~~~~~~~~~~~~~~~~~~~~~~~~~\pi-2\theta_{_{n}}=0\nonumber\\
\textnormal{or}~~~~~~~~~~~~~~~~~~\sin^4\theta_{_{n}}-\sin^2\theta_{_{n}}+(\pi/2-\theta_{_{n}})^2=0
\end{eqnarray}
Clearly, the only solution to these equations is
\begin{equation}
\theta_{_{n}}=\frac{\pi}{2}
\end{equation}
Therefore,
 \begin{eqnarray}
\frac{\pi}{2}=(\frac{1}{2}- \frac{1}{n})\pi\nonumber\\
\Rightarrow n\rightarrow\infty
\end{eqnarray}
The acceptable values for $n$ are $3\leq n<\infty$, this implies
$\pi/6\leq \theta_{_{n}}\leq \pi/2$. $d{\mathcal
A}_{_{n}}/d\theta_{_{n}}<0$ for $\theta_{_{n}}$ in [$\pi/3$,
$\pi/2$) and $d{\mathcal A}_{_{n}}/d\theta_{_{n}}> 0$ for
$\theta_{_{n}}$ in ($\pi/2$, $\infty$). This implies $A_{_{n}}$
attains its maximum at $\theta_{_{n}}= \pi/2$.

Next, we calculate the corresponding maximum area,
\begin{eqnarray}
{\mathcal
A}_{_{n-max}}=\frac{P^2}{8\pi}(\pi-2\theta_{_{n}})\tan{\theta_{_{n}}}|_{\theta_{_{n}}=\frac{\pi}{2}}\nonumber
\end{eqnarray}
This is in the indeterminate form: $0\times\infty$. We put the
above expression in the form $\frac{\infty}{\infty}$ and apply
L'Hopital's rule successively  to carry out the limiting
procedure:
\begin{eqnarray}
 {\mathcal
A}_{_{n-max}}=\frac{P^2}{8\pi}\lim _{{\theta_{_{n}}\rightarrow
\frac{\pi}{2}}}
\frac{\tan{\theta_{_{n}}}}{\frac{1}{\pi-2\theta_{_{n}}}}~~~~~~~\nonumber\\
=\frac{P^2}{16\pi}\lim _{{\theta_{_{n}}\rightarrow \frac{\pi}{2}}}
\frac{(\pi-2\theta_{_{n}})^2}{\cos^2\theta_{_{n}}}\nonumber\\
=\frac{P^2}{4\pi}\lim _{{\theta_{_{n}}\rightarrow \frac{\pi}{2}}}
\frac{(\pi-2\theta_{_{n}})}{\sin 2\theta_{_{n}}}~~~\nonumber\\
=-\frac{P^2}{4\pi}\lim _{{\theta_{_{n}}\rightarrow
\frac{\pi}{2}}}\frac{1}{\cos 2\theta_{_{n}}}~~~~\nonumber\\
=\frac{P^2}{4\pi}~~~~~~~~~~~~~~~~~~~~~~~\nonumber\\
=\pi r^2~~~(P=C=2\pi r)
\end{eqnarray}
Obviously the minimum area occurs for the case of a triangle:
\begin{eqnarray}
{\mathcal
A}_{_{n-min}}=\frac{P^2}{4n}\tan{\left[\left(\frac{1}{2}-
\frac{1}{n}\right)\pi\right]}_{n=3}\nonumber\\
=\frac{P^2}{12\sqrt 3}~~~~~~~~~~~~~~~~~~~~~~~
\end{eqnarray}

\section{Appendix B}
The following is the computer code that looks for the $\{n_{_1},
n_{_2},...,n_{_k}\}$ set of polygonal sides which maximizes the
value of total area in Eq.(78). This appendix also contains four
outputs of the program.
\newpage

\section{References}
1.~~~Varberg, D. and Purcell, E. J., \emph{Calculus with Analytic
Geometry}. 6th ed. New $~~~~~~$Jersey:Prentice Hall, 1992. \\
2.~~~Embry, M. R., Schell, J. F. and Thomas, J. P., \emph{Calculus
and Linear Algebra}. $~~~~~~$Philadelphia:W.B. Saunders Company, 1972.\\
3.~~~Apostol, T. M., \emph{Mathematical Analysis}.
Reading:Addison-Wesley, 1958.

\end{document}